\documentclass[11pt]{article}
\usepackage{graphicx}
\usepackage{amsfonts}
\usepackage{bbm}
\usepackage{mathrsfs}
 \usepackage{mathrsfs,amsmath,amssymb,amsthm}
\usepackage[dvips]{color}
 \usepackage{amssymb}
 \usepackage{amsbsy}
 \usepackage[dvips]{color}
\allowdisplaybreaks
 \theoremstyle{plain}
\theoremstyle{remark}  \newtheorem{remark}{\noindent\mbox{Remark}}
 \theoremstyle{plain}
 \theoremstyle{plain}\newtheorem{lemma}{\noindent\mbox{Lemma}}
\theoremstyle{plain} \newtheorem{theorem}{\noindent\mbox{Theorem}}
 \theoremstyle{plain}\newtheorem{proposition}{\noindent\mbox{Proposition}}
 \theoremstyle{plain}
\theoremstyle{definition} \newtheorem{definition}{\noindent\mbox{Definition}}

 \def\bq{\begin{equation}}
 \def\eq{\end{equation}}
 \def\eqn{\end{eqnarray}}
 \def\bqn{\begin{eqnarray}}
 \def\proof{\noindent{\it Proof.~~}}
 \def\qed{\hfill$\Box$\medskip}
 \def\rto{\rightarrow\infty}
 \def\z{\left}
 \def\y{\right}
 \def\no{\nonumber}

\begin{document}
 \title{\textbf{On the number of points skipped by a transient (1,2) random walk on the line}\footnote{Supported by National
Nature Science Foundation of China (Grant No. 11501008),  Nature Science Foundation of Anhui Province (Grant No. 1508085QA12) and Nature Science Foundation of Anhui Educational Committee (Grant No. KJ2014A085).}}                  

\author{  Hua-Ming \uppercase{Wang} \\
     {Department of Statistics, Anhui Normal University, Wuhu 241003, China}\\
    E-mail\,$:$ hmking@mail.ahnu.edu.cn}
\date{}
\maketitle%

\vspace{-.5cm}

\begin{center}
\begin{minipage}[c]{12cm}
\begin{center}\textbf{Abstract}\quad \end{center}

Consider a transient near-critical (1,2) random walk on the positive half line.  We give a criteria for the finiteness of the number of the skipped points (the points never visited) by the random walk. This result generalizes (partially) the criteria for the finiteness of the number of cutpoingts of the nearest neighbor random walk on the line by Cs\'aki, F\"olders, R\'ev\'esz [J Theor Probab (2010) 23: 624-638].
\vspace{0.2cm}

\vspace{0.2cm}

\end{minipage}
\end{center}

\section{Introduction}
Consider a random walk on the lattice of positive half line. Precisely, let $X=(X_k)_{k\ge0}$ be a Markov chain on $\mathbb Z^+=\{0,1,2,...\}$ whose transition matrix is give by
\begin{equation*}\label{rw}\begin{split}
  &P(X_{k+1}=n+2|X_k=n)=p_n\in (0,1),\\
&P(X_{k+1}=n-1|X_k=n)=q_n=1-p_n,\ n\ge1,\\
&P(X_{k+1}=2|X_k=0)=p_0=1.
\end{split}\end{equation*}
It is easy to see that if $p_n\equiv \frac{1}{3}$ for all $n\ge1,$ then $X$ is recurrent. We aim to study the near-critical case, that is, perturbing $p_n$ by letting $$p_n=\frac{1}{3}+r_n,\text{ with } r_n\in(-1/3,2/3).$$

Intuitionally, if the chain $X$ is recurrent, then every point will be visited at least once, since the left-oriented jumps are nearest neighbor. But if $X$ is transient, since the right-oriented jumps are always with size 2, $X$ may skip some points. The question is how many points will be skipped.

The finiteness of the number of skipped points is very sensitive to the decay speed of the perturbation $r_n.$ Our main purpose is to give a criteria to tell, under which conditions, the number of skipped points is finite or infinite.

Our motivation originates from  the nearest neighbor random walk studied in \cite{cfr,cfrb}. In \cite{cfrb},  the authors gave a criteria for the finiteness of the number of cutpoints of the transient nearest neighbor  random walk. For  the nearest neighbor setting, if we denote by $E_i$ the probability from $i$ to $i+1,$ then $\rho_i:=\frac{1-E_i}{E_i}$ is crucial for deriving everything. But for (1,2) random walk, the quantity $\frac{q_n}{p_n}$ itself does not work directly, instead, we need to use the $n$-th tail of a continued fraction, constructed from $\frac{q_n}{p_n},$ which is related to the representation of Markov chain by  circle and weights.

For $n\ge1,$ denote by $a_n=\frac{q_n}{p_n}$  and define \begin{equation}\label{ksn}\xi_n=\frac{1}{a_n}\z(1+\frac{1}{\xi_{n+1}}\y).\end{equation}
The solution $(\xi_n)_{n\ge1}$ of (\ref{ksn}) is not necessarily unique, but it is unique  if $\sum_{n=1}^\infty\frac{1}{a_n}=\infty,$ see \cite{dey98}, Lemma 1. We refer the reader to \cite{dey98}, page 204 for the following recurrence criteria of $X.$

 {\noindent\bf Recurrence criteria:} {\it the chain $(X_k)$ is transient if and only if $\sum_{n=1}^\infty\frac{1}{\xi_1\cdots \xi_n}<\infty$
  where $\xi_n, n\ge1$ is a solution of (\ref{ksn}) with $\xi_1>0.$}

  We now define the so-called skipped points which we concern. Throughout, we use the notation $\#\{\ \}$ to denote the number of elements in set $\{\ \}.$

  \begin{definition} If $\#\{n\ge0:X_n=k\}=0,$ we call the site $k$ a skipped points of $X.$
\end{definition}
\begin{remark}
  We remark that the skipped points in our setting are more or less similar to the cutpoints studied in \cite{cfrb}.  We call such a point $k$ a skipped point because it was never visited by the path of the walk
  $X.$
\end{remark}

If $p_n=1/3+r_n$ with $\lim_{n\rto}r_n=0,$ then $\sum_{n=1}^{\infty} a_n^{-1}=\sum_{n=1}^\infty p_n/q_n=\infty$ and hence  (\ref{ksn}) has a unique solution $\xi_n,n\ge 1.$  In this case, we set $$U_i=\frac{1}{\xi_i},\ i\ge 1.$$
For $n>m>0,$ introduce notation
\begin{equation*}
  D(m,n)=\sum_{j=m+1}^{n}\prod_{i=m+1}^{j-1}U_i,
\end{equation*} where the convention being that the empty product equals 0.
Let \begin{equation}\label{dmi}
  D(m):=\lim_{n\rto}D(m,n).
\end{equation}
We are now ready to state the main result.
\begin{theorem}\label{main}

   Suppose that $p_n=\frac{1}{3}+r_n,$ with $0\le r_n<\frac{2}{3},$ $r_n\downarrow 0$ and $r_n-r_{n+1}=O(r_n^2)$ as $n\rto.$ Let $\xi_n,n\ge1$ be  the solution of (\ref{ksn}), $U_i=\frac{1}{\xi_i},$  and $D(i)$ be the one defined in (\ref{dmi}). If $$\sum_{n=2}^\infty\frac{1}{D(n)\log n}<\infty,$$ then the Markov chain $X$ has at most finitely many skipped points almost surely.

   If there exists some $\delta>0$ such that $D(n)\le \delta n\log n$ for $n$ large enough and $$\sum_{n=2}^\infty\frac{1}{D(n)\log n}=\infty,$$
 then with a positive probability $p,$ $p\ge \frac{11}{27},$ the Markov chain $X$ has infinitely many skipped points.
      \end{theorem}

Similar to \cite{cfrb}, we give a criteria for the finiteness of the number of skipped points in terms of the perturbation $r_i.$
\begin{theorem}\label{cs} For $1\le i\le 3,$ let $r_i=\frac{1}{3}$ and for $i\ge4,$
  set $$r_i=\frac{1}{3}\z(\frac{1}{i}+\frac{1}{i(\log\log i)^\beta}\y).$$ Then
  if $\beta>1,$ almost surely $X$ has at most finitely many skipped points;
  if $\beta\le 1,$   with a positive probability $p,$ $p\ge \frac{11}{27},$  $X$ has infinitely many skipped points.
\end{theorem}

Theorem \ref{main} and  Theorem \ref{cs} are similar to those in \cite{cfrb}. We adopt the approach used in \cite{cfrb} to give the proofs. The key step is to calculate the probability that a point is a skipped point, because in our setting, the right-oriented jumps are not nearest neighbor. The main idea is to split $\mathbb Z^+$ into nonintersecting layers, that is, let $\mathbb Z^+=\cup_{k=0}^{\infty}L_k$ with $L_k=\{2k,2k+1\},k=0,1,2,...$ Although it is involved, it is possible to calculate the probability that $L_k$ (or both $L_k$ and $L_j$) containing  a skipped points.

We arrange the paper as follows. We devote Section \ref{pre} to giving some preliminary results for $\xi_n,$  $D(n)$ and estimating the escaping probability of the random walk. In Section \ref{skc}, the probability of a layer $L_k$ (or both $L_k$ and $L_j$) containing a skipped point is  studied. At last, in Section \ref{pr}, we finish the proofs of the main results.

\section{Preliminary results} \label{pre}
For $n\ge 1$ set $p_n=1/3+r_n$ with $0<r_n<2/3$ and $\lim_{n\rto}r_n=0.$ Then we have the following facts:

\noindent (a) $a_n=\frac{q_n}{p_n}=\frac{2/3-r_n}{1/3+r_n}=2-9r_n+O(r_n^2);$

\noindent (b) $\lim_{n\rto}a_n=2$ and thus (\ref{ksn}) has a unique solution $(\xi_n)_{n\ge1}.$

\noindent Throughout, the ``large O notation" $O(g(i))$ means $\lim_{i\rto}O(g(i))/g(i)=C$ for some constant $C,$ and the ``little o notation" $o(g(i))$ means $\lim_{i\rto}o(g(i))/g(i)=0.$

{\it \subsection{  Asymptotics of the continued fraction $\xi_{n},\ n\ge1.$}}

Iterating (\ref{ksn}), we have
\begin{equation}\label{ksr}
  \xi_n=\frac{1}{a_n}\z(1+\dfrac{a_{n+1}}{1+\dfrac{a_{n+2}}{1+...}}\y).
\end{equation}
Among the traditional notations of continued fractions,
\begin{equation}\label{cfk}
\mathrm K_{n=1}^\infty(a_n|1):=\dfrac{a_1}{1+\dfrac{a_2}{{1+\dots}}}
\end{equation}
is used to denote the continued fraction
and $$f^{(n)}=\mathrm K_{m=n+1}^\infty(a_k|1):=\dfrac{a_{n+1}}{1+\dfrac{a_{n+2}}{{1+\dots}}}$$ denotes the $n$th tail of $\mathrm K_{n=1}^\infty(a_n|1).$
 We can now write  (\ref{ksr}) as
\begin{equation}\label{ksnb}
  \xi_n=\frac{1}{a_n}\z(1+f^{(n)}\y).
\end{equation}
It is easy to see that \begin{equation}\label{fnp}
  f^{(n)}=\frac{a_{n+1}}{1+f^{(n+1)}}.
\end{equation}
By the Seidel-Stern Theorem(\cite{lw92}, page 98), the continued fraction in (\ref{cfk}) converges. The fact $\lim_{n\rto}a_n=2$ implies that (\cite{lw92}, page 151) \begin{equation}\label{fc}
  \lim_{n\rto}f^{(n)}=f:=\mathrm K_{n=1}^\infty(2|1)=1.
\end{equation}
\begin{lemma}\label{xie}
Suppose that $r_n\downarrow 0$ and $r_n-r_{n+1}=O(r_n^2)$ as $n\rto.$
Let $(\xi_n)_{n\ge1}$ be the unique solution of (\ref{ksn}). Then
  \begin{equation}\label{xir}\xi_{n}^{-1}=1-3r_{n}+O(r_{n}^2).\end{equation}
  Moreover, for some $n_0>0,$ $\xi_n^{-1}, n\ge n_0$ is increasing in $n.$
\end{lemma}
\proof In view of (\ref{fc}) we can assume that $$f^{(n)}=1+br_{n+1}+O(r_{n+1}^2).$$ Recall that $a_n=2-9r_n+O(r_n^2).$ Since $r_n-r_{n+1}=O(r_n^2),$ it follows from  (\ref{fnp}) that \begin{align*}
 f^{(n)}(1+f^{(n+1)})&=2+b(2r_{n+1}+r_{n+2})+O(r_{n+1}^2) \\
 &=2+3br_{n+1} +O(r_{n+1}^2)=2-9r_{n+1}+O(r_{n+1}^2).
\end{align*}
Consequently, we have $b=-3$ and thus  $$f^{(n)}=1-3r_{n+1}+O(r_{n+1}^2).$$
Comparing (\ref{ksn}) and (\ref{ksnb}), we conclude that
$$\xi_{n+1}^{-1}=f^{(n)}=1-3r_{n+1}+O(r_{n+1}^2)$$
which proves (\ref{xir}). Moreover, with $C$ some constant, we  have by (\ref{xir}) that
$\xi_{n}^{-1}=1-3r_{n}+Cr_{n}^2.$ Since $r_n\downarrow0,$ then there exists an $n_0>0$ such that for  $n>n_0,$ $\xi_n^{-1}-\xi_{n+1}^{-1}>0.$
  \qed

\subsection{Escaping probability}

For $i\ge 1,$ write $U_i=\frac{1}{\xi_i}.$ Then by Lemma \ref{xie}, we have $$U_i=1-3r_i+O(r_i^2)=e^{-3r_i+O(r_i^2)}.$$
Recall that for $n>m>0,$
$D(m,n)=\sum_{j=m+1}^{n}\prod_{i=m+1}^{j-1}U_i,$ and for $m\ge 1,$ $$D(m):=\lim_{n\rto}D(m,n).$$
It is easy to see that $D(n)$ is increasing in $n$ whenever $\xi_n$ is decreasing in $n.$  Also we have
$D(n)=1+U_{n+1}D_{n+1}$ and for $n\ge m+2,$ $$D(m,n)=D(m)\Big(1-\prod_{i=m}^{n-1}\Big(1-\frac{1}{D(i)}\Big)\Big).$$
To sum up, we have the following lemma.
\begin{lemma}\label{de} Suppose that $r_n\downarrow 0$ and $r_n-r_{n+1}=O(r_n^2)$ as $n\rto.$  Then with $n_0$ the one in Lemma \ref{xie}, $D(n),\ n\ge n_0$ is increasing in $n$ and
for fixed $n>m,$ $\frac{D(m,n)}{D(m)}$ is decreasing in $m.$

\end{lemma}
For integers $1\le a\le b\le c,$ let
\begin{equation*}
  P(a,b,c)=P(X \text{ hits }[0,a] \text{ before }[c,\infty]\big|X_0=b).
\end{equation*}
We have the following estimations of the escaping probabilities.
\begin{lemma}\label{ecpb} For any integers $1\le a\le b\le c,$
\begin{equation*}\frac{\sum_{i=b}^{c-1}U_{a+1}\cdots U_i}{1+\sum_{i=a+1}^{c-1}U_{a+1}\cdots U_i}\le P(a,b,c)\le \frac{\sum_{i=b}^{c}U_{a+1}\cdots U_i}{1+\sum_{i=a+1}^{c}U_{a+1}\cdots U_i}.\end{equation*}
\end{lemma}
The lemma can be proved by a space reversal argument of proof of the  corresponding lemma for (2,1) random walk given in Letchikov \cite{leta} (see Lemma 1, page 230 therein).

By Lemma \ref{ecpb}, we have the following facts:
\begin{align}
  &\frac{1+U_{a+1}}{D(a,c+1)}\le 1-P(a,a+2,c)\le \frac{1+U_{a+1}}{D(a,c)},\ c> a+2;\label{epa}\\
  & \frac{1}{D(a,c+1)}\le 1-P(a,a+1,c)\le \frac{1}{D(a,c)},\ c> a+1;\label{epb}\\
  &1-P(a,a+1,\infty)=\frac{1}{D(a)};\label{epc}\\
  &1-P(a,a+2,\infty)=\frac{1+U_{a+1}}{D(a)}\label{epd}.
\end{align}


\section{Skipped points} \label{skc}

Set $L_k=\{2k,2k+1\}, \ k\ge 0.$ Then $\mathbb Z_+=\bigcup_{k\ge0} L_k.$
Denote by
\begin{align*}
  C^S=\{k\ge1: L_k \text{ contains a skipped point}\}.
\end{align*}
\begin{proposition}\label{estc} Suppose that $r_n\downarrow 0$ and $r_n-r_{n+1}=O(r_n^2)$ as $n\rto.$ Then for any $\epsilon>0,$ there exists some $k_0>0$ such that \text{for }$k\ge k_0,$
\begin{align}
     &\frac{p_{2k}}{D(2k+1)}\le P(k\in C^S)\le \frac{1}{D(2k)} ,\label{skpa}
     \end{align}
     and for $k>j\ge k_0,$
     \begin{align}
     P(j\in& C^S, k\in C^S) \ge  (1+\epsilon)^{-1}p_{2j}p_{2k}\frac{1}{D(2j,2k)}  \frac{1}{D(2k+1)},\label{skpl}\\
      P(j\in& C^S, k\in C^S)  \le \frac{27}{11} (1+\epsilon)^2P(j\in C^S)P(k\in C^S)\frac{D(2j+1)}{D(2j+1,2k)}.\label{skpb}
\end{align}
\end{proposition}
\proof For $k\ge 1,$ define $T_k=\inf\{n\ge 0:X_n\in L_k\},$ the time  the walk entering into $L_k=\{2k,2k+1\}$ for the first time. We denote by \begin{align*}
  &h_k(1)=P(X_{T_k}=2k),\\
  &h_k(2)=P(X_{T_k}=2k+1),\ k\ge1;\\
  &\eta_{k,j}(1)=P(X \text{ enters }[j+1,\infty) \text{ at }j+1 |X_0=k),\\
  &\eta_{k,j}(2)=P(X \text{ enters }[j+1,\infty) \text{ at }j+2 |X_0=k),\ 1\le k\le j.
\end{align*}
Write also
  \begin{align*}
 &Q(a,b,c)=P(X\text{ hits }[c,\infty) \text{ before } a |X_0=b),\\
    &Q_1(a,b,c)=P(X\text{ hits }[c,\infty) \text{ before } a \text{ at }c|X_0=b),\\
    &Q_2(a,b,c)=P(X\text{ hits }[c,\infty) \text{ before } a \text{ at }c+1|X_0=b).
  \end{align*}
  It is clear that \begin{equation}\label{pab}Q(a,b,c)=Q_1(a,b,c)+Q_2(a,b,c)=1-P(a,b,c).\end{equation}

  As a beginning, we prove (\ref{skpa}).
  Note that the event $\z\{k\in C^S\y\}$ occurs if and only if $L_k$ contains a skipped point. It is easy to see that
   on the event $\{X_{T_k}=2k\},$ $2k$ can not be a skipped point and that on the event $\{X_{T_k}=2k+1\},$ $2k+1$ can not be a skipped point. Therefore if we write \begin{align*}
    B_{k,1}&:=\{X_{T_k}=2k,2k+1\text{ is a skipped point} \},\\
    B_{k,2}&:= \{X_{T_k}=2k+1,2k\text{ is a skipped point} \},
  \end{align*}then  $\z\{k\in C^S\y\}=B_{k,1}\cup B_{k,2}.$
Thus using the strong Markov property, we have
  \begin{align}\label{kcs}
    P(k\in C^S)&=P(B_{k,1})+P(B_{k,2})\no\\
    &=h_k(1)\eta_{2k,2k}(2)(1-P(2k+1,2k+2,\infty))\no\\
    &\quad\quad+h_k(2)(1-P(2k,2k+1,\infty)).
    \end{align}
    For an upper bound, abandoning the term $\eta_{2k,2k}(2)$, since by Lemma \ref{de}, $D(n), n\ge n_0$ is increasing in $n,$  then for $k>k_1\equiv \frac{n_0}{2}$ we have
           \begin{align*}
     P(k\in C^S)\le h_k(1)\frac{1}{D(2k+1)}+h_k(2)\frac{1}{D(2k)}\le \frac{1}{D(2k)}.
  \end{align*}
  For a lower bound, it follows from (\ref{kcs}) that for $k\ge k_1,$
   \begin{align}\label{kcsl}
    P(k\in C^S)& \ge [h_k(1)\eta_{2k,2k}(2)+h_k(2)]\frac{1}{D(2k+1)}.
    \end{align}
Since $\eta_{2k,2k}(2)\ge p_{2k},$ the lower bound in (\ref{skpa}) follows.

Next, we proceed to prove (\ref{skpl}) and (\ref{skpb}). Note that at most one site in $L_k$ is a skipped point, since at least one point in $L_k$ should be visited.
Write \begin{align*}
  &E_{jk}^{(1)}=\{\text{both }2j+1\text{ and }2k+1\text{ are skipped points}\},\\
    &E_{jk}^{(2)}=\{\text{both }2j+1\text{ and }2k\text{ are skipped points}\},\\
   &E_{jk}^{(3)}=\{\text{both }2j\text{ and }2k\text{ are skipped points}\},\\
        &E_{jk}^{(4)}=\{\text{both }2j\text{ and }2k+1\text{ are skipped points}\}.
\end{align*}
We calculate the probability of the event $E_{jk}^{(1)}$ in detail.
Indeed $E_{jk}^{(1)}$ occurs if and only if all of the following events occur:
\begin{center}
 \begin{minipage}[c]{11cm}
   {\bf (i)} the walk hits $L_j$ at $2j;$
   {\bf (ii)} starting from $2j,$ it hits $[2j+1,\infty)$ at $2j+2;$
      {\bf(iii)} starting from $2j+2,$ it hits $L_k$ at $2k$ before it returns to $2j+1;$
      {\bf(iv)} starting from $2k$ it hits $[2k+1,\infty)$ at $2k+2$ before it returns to $2j+1;$
      {\bf(v)} starting from $2k+2,$ it will never return to $2k+1.$
 \end{minipage}
\end{center}
Hence we have
\begin{align*}
  P\z(E_{jk}^{(1)}\y)=h_j(1)\eta_{2j,2j}(2)&Q_1(2j+1,2j+2,2k)\\
\times& Q_2(2j+1,2k,2k+1)(1-P(2k+1,2k+2,\infty)).
\end{align*}
 The probabilities of $E_{jk}^{(i)}, i=2,3,4 $ can be calculated similarly.
 We have
 \begin{align}
   P(j&\in C^S, k\in C^S)=\sum_{i=1}^4 P\z(E_{jk}^{(i)}\y)\no\\
   &=h_j(1){\eta_{2j,2j}(2)}Q_1(2j+1,2j+2,2k)\no\\
&\quad\quad\quad\quad\quad\quad\times {Q_2(2j+1,2k,2k+1)}(1-P(2k+1,2k+2,\infty))\no\\
&\quad\quad+h_j(1){\eta_{2j,2j}(2)}Q_2(2j+1,2j+2,2k)(1-P(2k,2k+1,\infty))\no\\
&\quad\quad+h_j(2)Q_2(2j,2j+1,2k)(1-P(2k,2k+1,\infty))\no\\
&\quad\quad+h_j(2)Q_1(2j,2j+1,2k){Q_2(2j,2k,2k+1)}\no\\
&\quad\quad\quad\quad\quad\quad\times(1-P(2k+1,2k+2,\infty)). \label{jki}
 \end{align}
By  Lemma \ref{de} and (\ref{epa})-(\ref{epd}), it follows that for $k> j >k_1,$
 \begin{align}\label{jke}
   P(j&\in C^S, k\in C^S)\no\\
   &\le h_j(1){\eta_{2j,2j}(2)}\frac{1}{D(2k)}\no\\
   &\quad\quad \times \{Q_1(2j+1,2j+2,2k){Q_2(2j+1,2k,2k+1)}+Q_2(2j+1,2j+2,2k)\}\no\\
&\quad\quad+ h_j(2)\frac{1}{D(2k)}\{Q_1(2j,2j+1,2k){Q_2(2j,2k,2k+1)}+Q_2(2j,2j+1,2k)\}\no\\
 &\le [h_j(1){\eta_{2j,2j}(2)}+h_j(2)]\frac{1}{D(2j+1)}[h_k(1){\eta_{2k,2k}(2)}+h_k(2)]\frac{1}{D(2k+1)}\no\\
 & \quad\quad\times {D(2j+1)}\frac{D(2k+1)}{D(2k)}\max\z\{Q(2j+1,2j+2,2k),Q(2j,2j+1,2k)\y\}\no\\
 &\quad\quad\times\max\z\{\frac{A(j,k)}{h_k(1){\eta_{2k,2k}(2)}+h_k(2)},\frac{B(j,k)}{h_k(1){\eta_{2k,2k}(2)}+h_k(2)}\y\}\no\\
 &\le P(j\in C^S)P(k\in C^S)\frac{D(2k+1)}{D(2k)}\frac{D(2j+1)}{D(2j)}\frac{D(2j+1)}{D(2j+1,2k)}\no\\
 &\quad\quad\times\max\z\{\frac{A(j,k)}{h_k(1){\eta_{2k,2k}(2)}+h_k(2)},\frac{B(j,k)}{h_k(1){\eta_{2k,2k}(2)}+h_k(2)}\y\}
 \end{align}
 where \begin{equation}\label{aij}A(j,k)=\frac{Q_1(2j+1,2j+2,2k)}{Q(2j+1,2j+2,2k)}{Q_2(2j+1,2k,2k+1)}+\frac{Q_2(2j+1,2j+2,2k)}{Q(2j+1,2j+2,2k)}\end{equation}
 and \begin{equation}\label{bij}B(j,k)=\frac{Q_1(2j,2j+1,2k)}{Q(2j,2j+1,2k)}{Q_2(2j,2k,2k+1)}+\frac{Q_2(2j,2j+1,2k)}{Q(2j,2j+1,2k)}.\end{equation}

It is easy to see that \begin{equation}\label{abe}
  A(j,k)\le 1,\ B(j,k)\le 1,
\end{equation} and  there exists a $k_2>0,$ such that for $k\ge k_2,$
\begin{equation}\label{etke}
  \eta_{2k,2k}(2) \ge \frac{11}{27}.
\end{equation}

Moreover, for any $\epsilon>0,$ since $\lim_{n\rto}D(n)=\infty$ (see \cite{cfr} Lemma 2.1, page 103) and $D(n)=1+U_{n+1}D(n+1),$  we can find $k_3>0$ such that for $k>j>k_3,$
\begin{equation}\label{dde}
  \frac{D(2k+1)}{D(2k)}\le 1+\epsilon,\ \frac{D(2j+1)}{D(2j)}\le 1+\epsilon.
\end{equation}
Consequently, substituting (\ref{abe})-(\ref{dde}) into the rightmost hand of (\ref{jke}), we conclude that for $k>j\ge k_0:=k_1\vee k_2\vee k_3,$
 \begin{align*}
   P(j\in C^S, k\in C^S)\le \frac{27}{11} (1+\epsilon)^2P(j\in C^S)P(k\in C^S)\frac{D(2j+1)}{D(2j+1,2k)}.
 \end{align*}

To prove (\ref{skpl}), by inserting some terms (those underlined) into (\ref{jki}), we have
 \begin{align*}
   P(j&\in C^S, k\in C^S)=\sum_{i=1}^4 P(E_{jk}^{(i)})\\
   &\ge h_j(1)\eta_{2j,2j}(2)Q_1(2j+1,2j+2,2k)\\
&\quad\quad\quad\quad\quad\quad\times Q_2(2j+1,2k,2k+1)(1-P(2k+1,2k+2,\infty))\\
&\quad\quad+h_j(1)\eta_{2j,2j}(2)Q_2(2j+1,2j+2,2k)\\ &\quad\quad\quad\quad\quad\quad\times\underline{Q_2(2j+1,2k,2k+1)}(1-P(2k,2k+1,\infty))\\
&\quad\quad+h_j(2)\underline{\eta_{2j,2j}(2)}Q_2(2j,2j+1,2k)\\
&\quad\quad\quad\quad\quad\quad\times \underline{Q_2(2j,2k,2k+1)}(1-P(2k,2k+1,\infty))\\
&\quad\quad+h_j(2)\underline{\eta_{2j,2j}(2)} Q_1(2j,2j+1,2k)\\
&\quad\quad\quad\quad\quad\quad\times Q_2(2j,2k,2k+1)(1-P(2k+1,2k+2,\infty)).
 \end{align*}
 It is easy to see that $Q_2(2j,2k,2k+1)\ge p_{2k},$ $\ Q_2(2j+1,2k,2k+1)\ge p_{2k}$ and $\eta_{2j,2j}(2)\ge p_{2j}.$ Therefore using Lemma \ref{de} and (\ref{dde}), for $k>j\ge k_0,$ we have
\begin{align*}
   P(j\in C^S, k\in C^S)&\ge p_{2j}p_{2k}\min\z\{\frac{1}{D(2j,2k)}, \frac{1}{D(2j+1,2k)}\y\}  \frac{1}{D(2k+1)}\\
   &\ge p_{2j}p_{2k}\frac{D(2j)}{D(2j+1)}\frac{1}{D(2j,2k)}  \frac{1}{D(2k+1)}\\
   &\ge (1+\epsilon)^{-1}p_{2j}p_{2k}\frac{1}{D(2j,2k)}  \frac{1}{D(2k+1)}.
 \end{align*}
  The proposition is proved.\qed

\section{Proofs of the main results} \label{pr}

\subsection{Proof of Theorem \ref{main}}

Basically speaking, the proof of Theorem \ref{main} follows the same line with that of \cite{cfrb}, Theorem 1.1. But there are some differences in details, so we repeat the proof with emphasis on the differences.
We now begin to prove Theorem \ref{main}. In the proof, except being otherwise stated, $c$ is a constant which may change line by line.

 We prove firstly the convergent part.  Define $$C_{j,k}=\{x:2^j<x\le 2^{k}, x\in C^S\}$$ and denote the cardinality of $C_{j,k}$  by $A_{j,k}\equiv|C_{j,k}|:=\#\{x:x\in C_{j,k}\}.$
     Write $$a_m:=P(A_{m,m+1}>0),\ b_m=\sum_{i=1}^{2^{m-1}}\min_{2^{m}< k \le 2^{m+1}} \frac{1}{D(2(k-i),2k+1)D(2k)}.$$  Let $l_m$ be the largest $k\in C_{m,m+1}$ if $C_{m,m+1}\neq \phi.$
      Then using (\ref{skpl}) we have for $m$  large enough,
       \begin{align*}
        \sum_{j=2^{m-1}+1}^{2^{m+1}}&P(j\in C^S)=E(A_{m-1,m+1})\\
              &\ge \sum_{k=2^{m}+1}^{2^{m+1}}P(A_{m,m+1}>0,l_m=k)E(A_{m-1,m+1}|A_{m,m+1}>0,l_m=k)\\
              &\ge\sum_{k=2^{m}+1}^{2^{m+1}}P(A_{m,m+1}>0,l_m=k)\sum_{i={2^{m-1}+1}}^k P(L_i\in C^S ,L_k\in C^S)\\
              &\ge cP(A_{m,m+1}>0)\min_{2^{m}< k \le 2^{m+1}} \sum_{i=1}^{2^{m-1}}\frac{1}{D(2(k-i),2k+1)D(2k)}\\
        &=c a_mb_m.
      \end{align*}
      It is easy to see that $U_m\le 1$ for large $m.$ So we have for $m$ large enough, $D(m,n)\le (n-m)$ and hence $$b_m\ge \sum_{i=1}^{2^{m-1}}\frac{1}{i+1}\ge cm.$$
      Consequently,  using Proposition \ref{estc}, we have
      \begin{align*}
        \sum_{m=1}^\infty &P(A_{m,m+1}>0)=\sum_{m=1}^\infty a_m\le\sum_{m=1}^\infty\frac{c}{b_m}\sum_{j=2^{m-1}+1}^{2^{m+1}}P(j\in C^S)\\
        &\le \sum_{m=1}^\infty\frac{c}{m}\sum_{j=2^{m-1}+1}^{2^{m+1}} \frac{1}{D(2j+1)}
          = \sum_{m=1}^\infty\frac{c}{m}\sum_{j=2^{m}+3}^{2^{m+2}+1} \frac{1}{D(j)}\\
          &\le c\sum_{m=1}^\infty \sum_{j=2^{m}+3}^{2^{m+2}+1} \frac{1}{D(j)\log j} \le \sum_{m=1}^\infty \frac{1}{D(n)\log n}<\infty.
      \end{align*}
         An application of the Borel-Cantelli lemma yields that with probability one, only finitely many of the events $\{A_{m,m+1}>0\}$ occur. We conclude that the Markov chain $X$ has at most finitely many skipped points almost surely.

         Next we prove the divergent part. Set $$m_k=[k\log k], A_k=\{m_k\in C^S\}.$$
         Here and throughout, $[x]$ denotes  integer part of $x.$ Our purpose is to show that \begin{equation*}\label{aio}
           P(A_k,k\ge 1 \text{ occur infinitely often})=1.
         \end{equation*}
         Now fix $\epsilon>0.$          By Lemma \ref{de} and Proposition  \ref{estc}  we can find $k_0$ such that for $k\ge k_0,$
          $$P(A_k)\ge \frac{c}{D(2m_k+1)}=\frac{c}{D(2[k\log k]+1)}\ge \frac{c}{D([2k\log 2k] )}$$
          and  for $l>k>k_0,$
          \begin{align}\label{lka}
           P(A_kA_l)&=P(m_k\in C^S, m_l\in C^S)\no\\
           & \le \frac{27}{11} (1+\epsilon)^2 P(m_k\in C^S)P(m_l\in C^S)\frac{D(2m_k+1)}{D(2m_k+1,2m_l)}\no\\
                   &\le\frac{27}{11} (1+\epsilon)^2 \z\{\frac{D(2m_k+1,2m_l)}{D(2m_k+1)}\y\}^{-1} P(A_k)P(A_l).
         \end{align}
         Thus (see \cite{cfrb}, page 630, Lemma 2.2) we have
         \begin{align*}
           \sum_{k\ge k_0}P(A_k)\ge \sum_{k\ge k_0} \frac{c}{D([2k\log 2k])} =\infty.
         \end{align*}
         Note that
         \begin{align*}
         \frac{D(2m_k+1,2m_l)}{D(2m_k+1)}&=1-\prod_{i=2m_k+1}^{2m_l-1}\Big(1-\frac{1}{D(i)}\Big)\ge 1-\exp\Bigg\{-\sum_{i=2m_k+1}^{2m_l-1}\frac{1}{D(i)}\Bigg\}.
         \end{align*}
      Define \begin{equation}\label{dl}
        l_1=\min\Bigg\{l\ge k:\sum_{i=2m_k+1}^{2m_l-1}\frac{1}{D(i)}\ge \log\frac{1+\epsilon}{\epsilon}\Bigg\}.
      \end{equation}
          For $l\ge l_1,$
     $\z\{\frac{D(2m_k+1,2m_l)}{D(2m_k+1)}\y\}^{-1}\le \frac{1}{1-\exp\z\{-\sum_{i=2m_k+1}^{2m_l-1}\frac{1}{D(i)}\y\}}\le 1+\epsilon.$
      Therefore, by (\ref{lka}) we have
      \begin{align}\label{po}
           P(A_kA_l)\le\frac{27}{11}(1+\epsilon)^3  P(A_l)P(A_l), \text{ for } l\ge l_1.
         \end{align}
         Next consider $k< l<l_1.$
Note that for $0\le u\le \log\frac{1+\epsilon}{\epsilon},$ we have
 $ 1-e^{-u}\ge cu.$
Since $D(n),\ n\ge n_0$ is increasing in $n,$  then by (\ref{lka}), we have
\begin{align*}
           P(A_kA_l) &\le\frac{27}{11} (1+\epsilon)^2\z\{\frac{D(2m_k+1,2m_l)}{D(2m_k+1)}\y\}^{-1} P(A_l)P(A_l)\\
           &\le \frac{c P(A_k)P(A_l)}{\sum_{i=2m_k+1}^{2m_l-1}\frac{1}{D(i)}} \le \frac{c P(A_k)P(A_l)D(2m_l)}{m_l-m_k}\le \frac{c P(A_k)}{l\log l-k\log k}.
         \end{align*}
Summing over $k+1$ to $l_1-1,$ we get
\begin{align}\label{kl}
   \sum_{l=k+1}^{l_1-1}&P(A_kA_l)\le cP(A_k) \sum_{l=k+1}^{l_1-1}\frac{1}{l\log l-k\log k}\no\\
   &\le cP(A_k)\frac{1}{\log k}\sum_{l=k+1}^{l_1-1}\frac{1}{l-k}\le  cP(A_k)\frac{\log l_1 }{\log k}.
\end{align}
We claim that \begin{equation}\label{klg}
  \frac{\log l_1 }{\log k}\le\gamma
\end{equation}
with some constant $\gamma$ depending only on $\epsilon.$
Indeed by (\ref{dl}), we know that
        \begin{equation*}
          \sum_{i=2m_k+1}^{2m_l-1}\frac{1}{D(i)}< \log\frac{1+\epsilon}{\epsilon}\text{ for }k<l<l_1.
        \end{equation*}
This implies that, for large $k,$ we have $l<k^\gamma$ with $\gamma>\z(\frac{1+\epsilon}{\epsilon}\y)^{\delta}+\epsilon.$ If we assume the contrary that $l\ge k^\gamma,$ then
\begin{align*}
   \sum_{i=2m_k+1}^{2m_l-1}\frac{1}{D(i)}&\ge \frac{1}{\delta}\sum_{i=2m_k+1}^{2m_l-1}\frac{1}{i\log i}\ge \frac{1}{\delta} \z[\log\log (2m_l)-\log\log(2m_k+1)\y]\\
   &\ge\frac{1}{\delta} \log(\gamma-\epsilon)\ge\log\frac{1+\epsilon}{\epsilon},
\end{align*}
a contradiction. Hence (\ref{klg}) is proved. Substituting (\ref{klg}) into (\ref{kl}), we have
\begin{align}\label{klc}
   \sum_{l=k+1}^{l_1-1}P(A_kA_l)\le  cP(A_k).
\end{align}
Taking (\ref{po}) and (\ref{klc}) together, we have
\begin{align}\label{splk}
\sum_{k=k_0}^{N}\sum_{l=k+1}^NP(A_kA_l) \le\sum_{k=k_0}^{N}\sum_{l=k+1}^N\frac{27}{11} (1+\epsilon)^3 P(A_l)P(A_l)+ c\sum_{k=k_0}^{N}P(A_k) \no
\end{align}
Consequently, writing $H(\epsilon)=\frac{27}{11} (1+\epsilon)^3,$ we have
\begin{align*}
  \alpha_H&:=\liminf_{N\rto}\frac{\sum_{k=k_0}^{N}\sum_{l=k+1}^NP(A_kA_l)-\sum_{k=k_0}^{N}\sum_{l=k+1}^NH P(A_l)P(A_l)}{\z[\sum_{k=k_0}^{N}P(A_k)\y]^2}\\
  &\le \lim_{N\rto}\frac{c}{\sum_{k=k_0}^{N}P(A_k)}=0.
\end{align*}
By a version of Borel-Cantelli lemma (see Petrov \cite{pe04}, page 235),
\begin{align*}
  P(A_k,&k\ge1 \text{ occur infinitely often} )\ge P(A_k,k\ge k_0 \text{ occur infinitely often} )\\
  &\ge  \frac{1}{H+2\alpha_H}=\z(\frac{27}{11} (1+\epsilon)^3+2\alpha_H\y)^{-1}\ge \frac{11}{27(1-\epsilon)^3}.
\end{align*}
Letting $\epsilon\rightarrow0,$ we conclude that
\begin{equation*}
   P(A_k,k\ge1 \text{ occur infinitely often} )\ge \frac{11}{27}.
\end{equation*}
The proof of the divergent part is completed.\qed

\subsection{Proof of Theorem \ref{cs}}
Recall that in Theorem \ref{cs}, for $1\le i\le 3,$ $r_i=\frac{1}{3}$ and for $i\ge 4,$  $$r_i=\frac{1}{3}\z(\frac{1}{i}+\frac{1}{i(\log\log i)^\beta}\y),  \text{ with } \beta>0.$$

\begin{lemma}\label{imif}
We have $r_i\downarrow 0$ as and $r_i-r_{i+1}=O\z(\frac{1}{i^2}\y),$ as $i\rto.$
\end{lemma}
\proof It suffices to show that
\begin{equation}\label{ii}\frac{1}{i[\log\log i]^{\beta}}-\frac{1}{(i+1)[\log\log(i+1)]^{\beta}}=o\Big(\frac{1}{i^2}\Big)\end{equation}
as $i\rto.$ To this end, note that
\begin{align*}
 & \frac{1}{i[\log\log i]^{\beta}}-\frac{1}{(i+1)[\log\log(i+1)]^{\beta}}\sim\frac{(i+1)[\log\log(i+1)]^{\beta}-i[\log\log i]^{\beta}}{i^2[\log\log i]^{2\beta}}
\end{align*}
as $i\rto.$ Thus (\ref{ii}) follows if we can show  \begin{equation}\label{bt}
  \lim_{i\rto}\frac{(i+1)[\log\log(i+1)]^{\beta}-i[\log\log i]^{\beta}}{[\log\log i]^{2\beta}}=0.
\end{equation}
Indeed, we have
\begin{align*}
  [\log\log(i+1)]^\beta=[\log\log i]^\beta\z(1+{\log\z(1+\frac{\log(1+1/i)}{\log i}\y)}/{\log\log i}\y)^\beta.
\end{align*}
Write $s(i)=\log\z(1+\frac{\log(1+1/i)}{\log i}\y)/\log\log i$ then $s(i)\sim \frac{1}{i\log i\log\log i}$ and $(1-s(i))^\beta-1\sim \beta s(i)$ as $i\rto.$ Consequently,
\begin{align*}
  &\frac{(i+1)[\log\log(i+1)]^{\beta}-i[\log\log i]^{\beta}}{[\log\log i]^{2\beta}}\sim\frac{(i+1)[(1+s(i))^\beta-1]}{[\log\log i]^{\beta}}\\
&\quad\quad\sim \frac{\beta}{\log i[\log\log i]^{\beta+1}}\rightarrow 0,
\end{align*}
as $i\rto.$ Thus (\ref{bt}) is proved. \qed

 By Lemma \ref{xie} and Lemma \ref{imif} we have   $$U_i\equiv \frac{1}{\xi_i}= 1-3r_i+O(r_i^2)=e^{-3r_i+O(r_i^2)}.$$
Then going  verbatim as that of \cite{cfrb}, Theorem 5.1, we have $$c_1k(\log\log k)^\beta\le D(k)\le c_1k(\log\log k)^\beta.$$
Consequently, Theorem \ref{cs} follows from Theorem \ref{main}. \qed

\begin{remark}
  We guess that the number $p$ in the divergent part of Theorem \ref{main} should be $1.$ We only get a low bound $\frac{11}{27}$ because in Proposition \ref{estc}, we show that for $k>j>k_0,$ \begin{align*}
      P(j\in& C^S, k\in C^S)  \le \frac{27}{11} (1+\epsilon)^2P(j\in C^S)P(k\in C^S)\frac{D(2j+1)}{D(2j+1,2k)}.
\end{align*}
The number $ \frac{27}{11} $  in this upper bound can be replaced by $1,$ if one can prove that the `maximum' term in the rightmost hand of (\ref{jke}) is less than $1,$ that is, $$\max\z\{\frac{A(j,k)}{h_k(1){\eta_{2k,2k}(2)}+h_k(2)},\frac{B(j,k)}{h_k(1){\eta_{2k,2k}(2)}+h_k(2)}\y\}\le 1.$$ As a consequence, one can get an almost sure result for the divergent part.
\end{remark}

\vspace{0.5cm}

\noindent{\large{\bf \Large Acknowledgements:}} The program was carried out during my visit to Yunshyong Chow in Academia Sinica. Thank the support (including the financial support) from both Professor Chow and Academia Sinica. Special thanks also extent to Professor Wenming Hong for introducing me the Lamperti random walk.


\begin{thebibliography}{99}
\addtolength{\itemsep}{-0.5em}
\bibitem{cfr} Cs\'aki, E., F\"oldes, A., R\'ev\'esz, P., Transient nearest neighbor random walk on the line, J. Theor. Probab., Vol. 22, 100-122, 2009.
    \bibitem{cfrb} Cs\'aki, E., F\"oldes, A., R\'ev\'esz, P., On the number of cutpoints of transient nearest neighbor random walk on the line, J. Theor. Probab., Vol. 23, 624-638, 2010.
\bibitem{dey98}Derriennic, Y., Random walks with jumps in random environments (Examples of cycle and weight representations), Probability Theory and
Mathematical Statistics(Proceedings of the  7th Vilnius Conference), pp 199-212, Utrecht, VSP Press, 1999.

   \bibitem{leta} Letchikov, A. V., A limit theorem for a  random walk in a random environment, Theory Probab. Appl. Vol. 33(2) pp 228-238, 1988.
\bibitem{lw92} Lorentzen, L. and Waadeland, H., Continued fractions with applications, North-Holland Publishing Co., Amsterdam, 1992.
    \bibitem{pe04}Petrov, V.V., A generalization of the Borel-Cantelli lemma, Statist. Probab. Lett., Vol. 67, 233-239, 2004.
 \end{thebibliography}
\end{document}